# Characterizations of Slant Ruled Surfaces in the Euclidean 3-space

**Mehmet Önder, Onur Kaya**

*Celal Bayar University, Faculty of Arts and Sciences, Department of Mathematics, Muradiye Campus, 45047 Muradiye, Manisa, Turkey.*

E-mails: mehmet.onder@cbu.edu.tr, onur.kaya@cbu.edu.tr

**Abstract**

In this study, we give the relationships between the conical curvatures of ruled surfaces drawn by the unit vectors of the ruling, central normal and central tangent of a regular ruled surface in the Euclidean 3-space. We obtain the differential equations characterizing slant ruled surfaces and if the reference ruled surface is a slant ruled surface, we give some conditions for the surfaces drawn by the central normal and the central tangent vectors to be slant ruled surfaces.

**MSC:** 53A25, 14J26

**Keywords:** Slant ruled surface; Frenet frame; conical curvature.

## 1. Introduction

In differential geometry of curves and surfaces, special curves and surfaces have an important role and more applications. Generally, special curves are such curves whose curvatures satisfy some special conditions. One of the well-known of such curves is the helix curve in the Euclidean 3-space $E^3$ which is defined by the property that the tangent line of the curve makes a constant angle with a fixed straight line called the axis of the general helix [2]. Therefore, a general helix can be equivalently defined as one whose tangent indicatrix is a planar curve.

Recently, another special curve similar to helix, called slant helix, has been defined by Izumiya and Takeuchi [3]. They defined a slant helix such as the normal lines of the curve make a constant angle with a fixed direction and they have given a characterization of slant helix in the Euclidean 3-space $E^3$. Moreover, slant helices have been studied by some mathematicians and new types of these curves have been introduced. Kula and Yaylı have investigated spherical images, the tangent indicatrix and the binormal indicatrix of a slant helix and they have obtained that the spherical images of a slant helix are helices lying on unit sphere [5]. Kula and et al have introduced some new results characterizing slant helices in $E^3$ [6]. Ali has studied the position vectors of slant helices in the Euclidean 3-space [1]. Monterde has shown that for a curve with constant curvature and non-constant torsion the principal normal vector of the curve makes a constant angle with a fixed constant direction, i.e., the curve is a slant helix [7]. Later, Önder and et al have defined a new type of slant helix called $B_2$-slant helix in the Euclidean 4-space $E^4$ by saying that the second binormal vector of a space curve makes a constant angle with a fixed direction in $E^4$ [10].

In the case of surfaces, ruled surfaces are a kind of special surfaces which are generated by a continuous movement of a line along a curve. Önder has considered the notion of "slant helix" for ruled surfaces and defined slant ruled surfaces in $E^3$ by the property that the vectors of the Frenet frame of the surface make constant angles with fixed directions and he has given characterizations for a regular ruled surface to be a slant ruled surface [9]. He has also shown that the striction curves of developable slant ruled surfaces are helices or slant helices. Later, Önder and Kaya have defined Darboux slant ruled surfaces in $E^3$ such as the Darboux vector of the ruled surface makes a constant angle with a fixed direction and they have given characterizations for a regular ruled surface to a Darboux slant ruled surface [8].

In this work, we give new characterizations for slant ruled surfaces in $E^3$.

## 2. Ruled Surfaces in the Euclidean 3-space $E^3$

In this section, we give a brief summary of the geometry of ruled surfaces in $E^3$.

A ruled surface $S$ is a special surface generated by a continuous moving of a line along a curve and has the parametrization

$$\vec{r}(u,v) = \vec{f}(u) + v\vec{q}(u), \qquad (2.1)$$

where $\vec{f} = \vec{f}(u)$ is a regular curve in $E^3$ defined on an open interval $I \subset \mathbb{R}$ and $\vec{q} = \vec{q}(u)$ is a unit direction vector of an oriented line in $E^3$. The curve $\vec{f} = \vec{f}(u)$ is called base curve or generating curve of the surface and various positions of the generating lines $\vec{q} = \vec{q}(u)$ are called rulings. In particular, if the direction of $\vec{q}$ is constant, then the ruled surface is said to be cylindrical, and non-cylindrical otherwise.

Let $\vec{m}$ be unit normal vector of ruled surface $S$. Then we have

$$\vec{m} = \frac{\vec{r}_u \times \vec{r}_v}{\|\vec{r}_u \times \vec{r}_v\|} = \frac{(\dot{\vec{f}} + v\dot{\vec{q}}) \times \vec{q}}{\sqrt{\left\langle \dot{\vec{f}} + v\dot{\vec{q}},\ \dot{\vec{f}} + v\dot{\vec{q}} \right\rangle - \left\langle \dot{\vec{f}}, \vec{q} \right\rangle^2}}, \qquad (2.2)$$

where "dot" shows the derivative with respect to $u$. If $v$ infinitely decreases, then along a ruling $u = u_1$, the unit normal $\vec{m}$ approaches a limiting direction. This direction is called the asymptotic normal (central tangent) direction and from (2.2) defined by

$$\vec{a} = \lim_{v \to \pm\infty} \vec{m}(u_1, v) = \frac{\dot{\vec{q}} \times \ddot{\vec{q}}}{\|\ddot{\vec{q}}\|}.$$

The point at which the unit normal of $S$ is perpendicular to $\vec{a}$ is called the striction point (or central point) $C$ and the set of striction points of all rulings is called striction curve of the surface.

The vector $\vec{h}$ defined by $\vec{h} = \vec{a} \times \vec{q}$ is called central normal vector. Then the orthonormal system $\{C; \vec{q}, \vec{h}, \vec{a}\}$ is called Frenet frame of the ruled surface $S$ where $C$ is the central point and $\vec{q}$, $\vec{h}$, $\vec{a}$ are unit vectors of ruling, central normal and central tangent, respectively.

The set of all bound vectors $\vec{q}(u)$ at the point $O$ constitutes a cone which is called *directing cone* of the ruled surface $S$. The end points of unit vectors $\vec{q}(u)$ drive a spherical curve $k_q$ on the unit sphere $S_1^2$ and this curve is called *spherical image* of ruled surface $S$, whose arc length is denoted by $s_q$.

For the Frenet formulae of the ruled surface $S$ and of its directing cone with respect to the arc length $s_q$ we have

$$\begin{bmatrix} d\vec{q}/ds_q \\ d\vec{h}/ds_q \\ d\vec{a}/ds_q \end{bmatrix} = \begin{bmatrix} 0 & 1 & 0 \\ -1 & 0 & \kappa_q \\ 0 & -\kappa_q & 0 \end{bmatrix} \begin{bmatrix} \vec{q} \\ \vec{h} \\ \vec{a} \end{bmatrix}, \quad (2.3)$$

where $\kappa_q(s_q) = \|d\vec{a}/ds_q\|$ is called the conical curvature of the directing cone (For details [4]).

**Definition 2.1 ([10])** Let $S_q$ be a regular ruled surface in $E^3$ given by the parametrization

$$\vec{r}(s, v) = \vec{c}(s) + v\vec{q}(s), \quad \|\vec{q}(s)\| = 1,$$

where $\vec{c}(s)$ is striction curve of $S$ and $s$ is arc length parameter of $\vec{c}(s)$. Let the Frenet frame of $S$ be $\{\vec{q}, \vec{h}, \vec{a}\}$. Then $S$ is called a $q$-slant ($h$-slant or $a$-slant, respectively) ruled surface if the ruling (the vector $\vec{h}$ or the vector $\vec{a}$, respectively) makes a constant angle with a fixed non-zero direction $\vec{u}$ in the space, i.e.,

$$\langle \vec{q}, \vec{u} \rangle = \cos\theta = constant; \quad \theta \neq \frac{\pi}{2}, \tag{2.4}$$

$$(\langle \vec{h}, \vec{u} \rangle = \cos\theta = constant; \quad \theta \neq \frac{\pi}{2} \text{ or } \langle \vec{a}, \vec{u} \rangle = \cos\theta = constant; \quad \theta \neq \frac{\pi}{2}, \text{ respectively}).$$

**Theorem 2.1 ([8])** *Let $S_q$ be a regular ruled surface in $E^3$ with Frenet frame $\{\vec{q}, \vec{h}, \vec{a}\}$ and conical curvature $\kappa_q \neq 0$. Then $S_q$ is an $h$-slant ruled surface if and only if the function*

$$\frac{\kappa'_q}{\left(1 + \kappa_q^2\right)^{3/2}}, \tag{2.5}$$

*is constant.*

**Theorem 2.2** *Let $S_q$ be a regular ruled surface in $E^3$ with Frenet frame $\{\vec{q}, \vec{h}, \vec{a}\}$ and conical curvature $\kappa_q \neq 0$. Then $S_q$ is a $q$-slant ruled surface if and only if the function $\kappa_q$ is constant.*

***Proof.*** Let $S_q$ be a $q$-slant ruled surface in $E^3$. Then denoting by $\vec{u}$ the unit vector of fixed direction and by $\theta$ the angle between $\vec{q}$ and $\vec{u}$, the following equality is satisfied

$$\langle \vec{q}, \vec{u} \rangle = \cos\theta = constant. \tag{2.6}$$

By taking the derivative of (2.6) with respect to $s_q$ gives $\langle \vec{h}, \vec{u} \rangle = 0$. Therefore, the vector $\vec{u}$ lies on the plane spanned by the vectors $\vec{q}$ and $\vec{a}$, i.e.,

$$\vec{u} = (\cos\theta)\vec{q} + (\sin\theta)\vec{a}. \tag{2.7}$$

By differentiating (2.7) with respect to $s_q$ and considering that the direction of $\vec{u}$ is constant it follows

$$0 = (\cos\theta - \kappa_q \sin\theta)\vec{h},$$

and since $\vec{h} \neq \vec{0}$, we have $\kappa_q = \cot\theta$ is constant.

Conversely, let $\kappa_q = \cot\theta$ be constant. We define,

$$\vec{u} = (\cos\theta)\vec{q} + (\sin\theta)\vec{a}.$$

Differentiating the last equation and using $\kappa_q = \cot\theta$ is constant we get $\vec{u}' = 0$, i.e., $\vec{u}$ is a constant vector. On the other hand, $\langle \vec{q}, \vec{u} \rangle = \cos\theta = constant$. Then we conclude that $S_q$ is a $q$-slant ruled surface.

## 3. Frenet Formulae of the Ruled Surfaces Generated by the Central Normal Vector and Central Tangent Vector

In this section, we give the Frenet formulae of the ruled surfaces generated by the central normal vector $\vec{h}$ and the central tangent vector $\vec{a}$ of the Frenet frame $\{\vec{q}, \vec{h}, \vec{a}\}$ of a regular ruled surface $S_q$. We show the ruled surfaces generated by $\vec{h}$ and $\vec{a}$ by $S_h$ and $S_a$, respectively; and their Frenet formulae and conical curvatures by $\{\vec{q}_h, \vec{h}_h, \vec{a}_h\}$, $\kappa_h$ and $\{\vec{q}_a, \vec{h}_a, \vec{a}_a\}$, $\kappa_a$, respectively.

**Theorem 3.1** *Let $S_q$ be a regular ruled surface in $E^3$ with the Frenet frame $\{\vec{q}, \vec{h}, \vec{a}\}$ and with non-zero conical curvature $\kappa_q$. Then the relationships between the conical curvatures of the surfaces $S_q$ and $S_h$ is given by*

$$\kappa_h = \frac{\kappa_q'}{\left(1 + \kappa_q^2\right)^{3/2}}.$$

*Then the Frenet formulae of $S_h$ is*

$$\begin{bmatrix} d\vec{q}_h / ds_h \\ d\vec{h}_h / ds_h \\ d\vec{a}_h / ds_h \end{bmatrix} = \begin{bmatrix} 0 & 1 & 0 \\ -1 & 0 & \kappa_h \\ 0 & -\kappa_h & 0 \end{bmatrix} \begin{bmatrix} \vec{q}_h \\ \vec{h}_h \\ \vec{a}_h \end{bmatrix}. \tag{3.1}$$

*Proof.* For the parametrization of $S_h$ we have

$$\vec{r}_h(s, v) = \vec{c}(s) + v\vec{h}(s), \quad \|\vec{h}(s)\| = 1,$$

where $\vec{c}(s)$ is striction curve of $S_q$ and $s$ is arc length parameter of $\vec{c}(s)$. Since the Frenet frame of $S_h$ is given by $\{\vec{q}_h, \vec{h}_h, \vec{a}_h\}$, we can write

$$\vec{q}_h = \vec{h},$$

and if we use the Frenet formulae given by (2.3), we get the central normal vector of $S_h$ as,

$$\vec{h}_h = \frac{d\vec{q}_h}{ds_q}\frac{ds_q}{ds_h} = (-\vec{q}+\kappa_q \vec{a})\frac{ds_q}{ds_h}, \qquad (3.2)$$

where $s_h$ is the arc length parameter of the spherical curve drawn by $\vec{h}$. Since $\vec{h}_h$ is a unit vector, from (3.2) we have

$$\frac{ds_q}{ds_h} = \frac{1}{\sqrt{1+\kappa_q^2}}.$$

Then (3.2) becomes

$$\vec{h}_h = \frac{1}{\sqrt{1+\kappa_q^2}}(-\vec{q}+\kappa_q \vec{a}). \qquad (3.3)$$

Since $\vec{q}_h = \vec{h}$, from (3.3), the central tangent vector is

$$\vec{a}_h = \vec{q}_h \times \vec{h}_h = \frac{1}{\sqrt{1+\kappa_q^2}}(\vec{a}+\kappa_q \vec{q}),$$

and the conical curvature of $S_h$ is

$$\kappa_h = \left\|\frac{d\vec{a}_h}{ds_h}\right\| = \left\|\frac{d\vec{a}_h}{ds_q}\frac{ds_q}{ds_h}\right\| = \frac{\kappa_q'}{\left(1+\kappa_q^2\right)^{3/2}}, \qquad (3.4)$$

where $\kappa_q' = \dfrac{d\kappa_q}{ds_q}$ and we have,

$$\begin{bmatrix} d\vec{q}_h/ds_h \\ d\vec{h}_h/ds_h \\ d\vec{a}_h/ds_h \end{bmatrix} = \begin{bmatrix} 0 & 1 & 0 \\ -1 & 0 & \kappa_h \\ 0 & -\kappa_h & 0 \end{bmatrix} \begin{bmatrix} \vec{q}_h \\ \vec{h}_h \\ \vec{a}_h \end{bmatrix}.$$

where $\kappa_h = \dfrac{\kappa_q'}{\left(1+\kappa_q^2\right)^{3/2}}$.

From Theorem 3.1 we can give the following theorem:

**Theorem 3.2** *The surface $S_q$ is an $h$-slant ruled surface if and only if the surface $S_h$ is a $q$-slant ruled surface.*

**Theorem 3.3** *Let $S_q$ be a regular ruled surface in $E^3$ with the Frenet frame $\{\vec{q}, \vec{h}, \vec{a}\}$ and with non-zero conical curvature $\kappa_q$. Then the relationships between the conical curvatures of the surfaces $S_q$ and $S_a$ is given by $\kappa_a = \dfrac{1}{\kappa_q}$. Then the Frenet formulae of $S_a$ is*

$$\begin{bmatrix} d\vec{q}_a / ds_a \\ d\vec{h}_a / ds_a \\ d\vec{a}_a / ds_a \end{bmatrix} = \begin{bmatrix} 0 & 1 & 0 \\ -1 & 0 & \kappa_a \\ 0 & -\kappa_a & 0 \end{bmatrix} \begin{bmatrix} \vec{q}_a \\ \vec{h}_a \\ \vec{a}_a \end{bmatrix}. \tag{3.5}$$

*Proof.* For the parametrization of the surface $S_a$ we have

$$\vec{r}_a(s, v) = \vec{c}(s) + v\vec{a}(s), \quad \|\vec{a}(s)\| = 1,$$

where $\vec{c}(s)$ is striction curve of $S_q$ and $s$ is arc length parameter of $\vec{c}(s)$. Since the Frenet frame of $S_a$ is given by $\{\vec{q}_a, \vec{h}_a, \vec{a}_a\}$, we can write,

$$\vec{q}_a = \vec{a},$$

and use the Frenet formulae given by (2.3), we obtain the central normal vector $\vec{h}_a$ as follows,

$$\vec{h}_a = \frac{d\vec{q}_a}{ds_q} \frac{ds_q}{ds_a} = -\kappa_q \vec{h}_q \frac{ds_q}{ds_a}, \tag{3.6}$$

where $s_a$ is the arc length parameter of the spherical curve drawn by $\vec{a}$. Since $\vec{h}_a$ is a unit vector, from (3.6) we have the followings,

$$\frac{ds_q}{ds_a} = \frac{1}{\kappa_q}, \quad \vec{h}_a = -\vec{h}.$$

Then the central tangent vector of the surface $S_a$ is

$$\vec{a}_a = \vec{q}_a \times \vec{h}_a = \vec{q},$$

and the conical curvature of $S_a$ is

$$\kappa_a = \left\| \frac{d\vec{a}_a}{ds_a} \right\| = \left\| \frac{d\vec{a}_a}{ds_q} \frac{ds_q}{ds_a} \right\| = \frac{1}{\kappa_q}.$$

Therefore we have

$$\begin{bmatrix} d\vec{q}_a / ds_a \\ d\vec{h}_a / ds_a \\ d\vec{a}_a / ds_a \end{bmatrix} = \begin{bmatrix} 0 & 1 & 0 \\ -1 & 0 & \kappa_a \\ 0 & -\kappa_a & 0 \end{bmatrix} \begin{bmatrix} \vec{q}_a \\ \vec{h}_a \\ \vec{a}_a \end{bmatrix},$$

where $\kappa_a = \dfrac{1}{\kappa_q}$.

From Theorem 3.3 we have the following theorem:

**Theorem 3.4** *The surface $S_q$ is a $q$-slant ruled surface if and only if the surface $S_a$ is a $q$-slant ruled surface.*

## 4. Differential Equation Characterizations of Slant Ruled Surfaces in $E^3$

In this section, we give the differential equations characterizing slant ruled surfaces.

**Theorem 4.1** *Let $S_q$ be a regular ruled surface in $E^3$ with the Frenet frame $\{\vec{q}, \vec{h}, \vec{a}\}$ and non-zero conical curvature $\kappa_q$. Then, $S_q$ is a $q$-slant ruled surface if and only if the vector $\vec{q}$ satisfies the following differential equation,*

$$\vec{q}''' + (1 + \kappa_q^2)\vec{q}' = 0, \tag{4.1}$$

*where $\vec{q}' = d\vec{q}/ds_q$, $\vec{q}''' = d^3\vec{q}/ds_q^3$.*

***Proof.*** Let $S_q$ be a $q$-slant ruled surface in $E^3$. From the Frenet formulae given in (2.3) we have,

$$\vec{q}' = \vec{h}.$$

If we take the derivative of the last equation twice with respect to $s_q$ we obtain

$$\vec{q}''' = -\vec{h} + \kappa_q' \vec{a} - \kappa_q^2 \vec{h}. \tag{4.2}$$

Since $S_q$ is a $q$-slant ruled surface from Theorem 2.2, $\kappa_q = $ constant, i.e., $\kappa_q' = 0$. Hence, from (4.2) it follows

$$\vec{q}''' + (1 + \kappa_q^2)\vec{q} = 0,$$

which is desired.

Conversely, let us assume that the equation (4.1) holds. From the Frenet formulae in (2.3) we have,

$$\vec{h}' = -\vec{q} + \kappa_q \vec{a}. \tag{4.3}$$

Since $\kappa_q \neq 0$, from (4.3) we can write

$$\vec{a} = \frac{1}{\kappa_q}(\vec{h}' + \vec{q}), \qquad (4.4)$$

and by taking the derivative of (4.4) we get

$$\vec{a}' = \frac{1}{\kappa_q}\left(\vec{h}'' + \vec{q}'\right) - \frac{\kappa_q'}{\kappa_q^2}\left(\vec{h}' + \vec{q}\right).$$

By using the Frenet formulae given in (2.3) and the equation (4.1), from the last equation we obtain

$$\frac{\kappa_q'}{\kappa_q}\vec{a} = 0,$$

where $\kappa_q' = 0$ which gives us that $\kappa_q = $ constant. Then, from Theorem 2.2, $S_q$ is a $q$-slant ruled surface.

**Theorem 4.2** *Let $S_q$ be a regular ruled surface in $E^3$ with the Frenet frame $\{\vec{q}, \vec{h}, \vec{a}\}$ and non-zero conical curvature $\kappa_q$. Then, $S_q$ is a $q$-slant ruled surface if and only if the vector $\vec{h}$ satisfies the following differential equation*

$$\vec{h}'' + (1 + \kappa_q^2)\vec{h} = 0. \qquad (4.5)$$

***Proof.*** Let $S_q$ be a $q$-slant ruled surface in $E^3$. From the Frenet formulae in (2.3) we have

$$\vec{h}' = -\vec{q} + \kappa_q \vec{a}.$$

By taking the derivative of the last equation with respect to $s_q$ we obtain

$$\vec{h}'' = -\vec{h} + \kappa_q' \vec{a} - \kappa_q^2 \vec{h}. \qquad (4.6)$$

Since $S_q$ is a $q$-slant ruled surface, from Theorem 2.2, $\kappa_q = $ constant. Hence, from (4.6) it follows,

$$\vec{h}'' + (1 + \kappa_q^2)\vec{h} = 0.$$

Conversely, let the equation (4.5) holds. From (2.3) we have

$$\vec{h}' = -\vec{q} + \kappa_q \vec{a}.$$

By taking the derivative of the last equation with respect to $s_q$ we get

$$\vec{q}' = -\vec{h}'' + \kappa_q' \vec{a} - \kappa_q^2 \vec{h},$$

and finally, by using (2.3) and (4.5) we obtain

$$\vec{h} = \vec{h} + \kappa_q' \vec{a},$$

which gives us $\kappa_q' = 0$, and so, $\kappa_q = $ constant. Then Theorem 2.2 gives that $S_q$ is a $q$-slant ruled surface in $E^3$.

**Theorem 4.3** *Let $S_q$ be a regular ruled surface in $E^3$ with the Frenet frame $\{\vec{q}, \vec{h}, \vec{a}\}$ and with non-zero conical curvature $\kappa_q$. Then, $S_q$ is a $q$-slant ruled surface if and only if the central tangent vector $\vec{a}$ satisfies the following differential equation*

$$\vec{a}''' + (1 + \kappa_q^2)\vec{a}' = 0. \tag{4.7}$$

***Proof.*** Since $S_q$ is a $q$-slant ruled surface, $\kappa_q = $ constant, i.e., $\kappa_q' = 0$. Then, from (2.3) we have

$$\vec{a}' = -\kappa_q \vec{h}. \tag{4.8}$$

If we take the derivative of the equation (4.8) and use the Frenet formulae given by (2.3) it follows

$$\vec{a}''' + (1 + \kappa_q^2)\vec{a}' = 0.$$

Conversely, let (4.7) holds. From the Frenet formulae in (2.3) we have

$$\vec{a}' = -\kappa_q \vec{h}.$$

If we differentiate the last equation twice, we get

$$2\kappa_q' \vec{q} - \kappa_q'' \vec{h} - 3\kappa_q \kappa_q' \vec{a} = 0. \tag{4.9}$$

Since the Frenet frame is linearly independent, from (4.9) we obtain the following system

$$2\kappa_q' = 0, \ \kappa_q'' = 0, \ 3\kappa_q \kappa_q' = 0,$$

which gives us that $\kappa_q = $ constant. Therefore, $S_q$ is a $q$-slant ruled surface.

In the following theorems, we give the differential equation characterizations of the surface $S_q$ by means of the Frenet vectors of ruled surfaces $S_h$ and $S_a$.

**Theorem 4.4** *Let $S_q$ be a regular ruled surface in $E^3$ with the Frenet frame $\{\vec{q}, \vec{h}, \vec{a}\}$ and with non-zero conical curvature $\kappa_q$. Then, $S_q$ is a $h$-slant ruled surface if and only if the ruling vector $\vec{q}_h$ of the ruled surface $S_h$ satisfies the following equation*

$$\frac{d^3\vec{q}_h}{ds_h^3}+(1+\kappa_h^2)\frac{d\vec{q}_h}{ds_h}=0 \tag{4.10}$$

**Proof.** Let $S_q$ be an $h$-slant ruled surface. Then from Theorem 2.1 we have that $\dfrac{\kappa'_q}{\left(1+\kappa_q^2\right)^{3/2}}$ is constant. From (3.1) we have

$$\frac{d\vec{q}_h}{ds_h}=\vec{h}_h. \tag{4.11}$$

Differentiating (4.11) twice and using that $\kappa_h=\dfrac{\kappa'_q}{\left(1+\kappa_q^2\right)^{3/2}}$ is constant, we get

$$\frac{d^3\vec{q}_h}{ds_h^3}+(1+\kappa_h^2)\frac{d\vec{q}_h}{ds_h}=0,$$

which is desired.

Conversely, let us assume that the equation (4.10) holds. From (3.1) we have

$$\frac{d\vec{q}_h}{ds_h}=\vec{h}_h.$$

By taking the derivative of the last equation with respect to $s_h$ and using (4.10) it follows $\kappa'_h\vec{a}_h=0$ which gives that $\kappa_h=\dfrac{\kappa'_q}{\left(1+\kappa_q^2\right)^{3/2}}$ is constant, and from Theorem 2.1 we have that $S_q$ is a $h$-slant ruled surface.

**Theorem 4.5** *Let $S_q$ be a regular ruled surface in $E^3$ with the Frenet frame $\{\vec{q},\vec{h},\vec{a}\}$ and with non-zero conical curvature $\kappa_q$. Then, $S_q$ is an $h$-slant ruled surface if and only if the central normal vector $\vec{h}_h$ of the ruled surface $S_h$ satisfies the equation*

$$\frac{d^2\vec{h}_h}{ds_h^2}+(1+\kappa_h^2)\vec{h}_h=0 \tag{4.12}$$

*where $\kappa_h$ is the conical curvature of the surface $S_h$.*

**Proof.** Let $S_q$ be an $h$-slant ruled surface. From (3.1) we have

$$\frac{d\vec{h}_h}{ds_h}=-\vec{q}_h+\kappa_h\vec{a}_h.$$

If we take the derivative of the last equation we obtain

$$\frac{d^2\vec{h}_h}{ds_h^2} = -\vec{h}_h + \frac{d\kappa_h}{ds_h}\vec{a}_h - \kappa_h^2\vec{h}_h. \qquad (4.13)$$

Since $S$ is an $h$-slant ruled surface then, from Theorem 2.1 and Theorem 3.1.

$$\kappa_h = \frac{\kappa'_q}{\left(1+\kappa_q^2\right)^{3/2}},$$

is constant. Then from (4.13) we have

$$\frac{d^2\vec{h}_h}{ds_h^2} + (1+\kappa_h^2)\vec{h}_h = 0. \qquad (4.14)$$

Conversely, let us assume that (4.12) holds. From (3.1) we have that

$$\frac{d^2\vec{h}_h}{ds_h^2} = -\vec{q}'_h + \frac{d\kappa_h}{ds_h}\vec{a}_h - \kappa_h^2\vec{h}_h. \qquad (4.15)$$

If we substitute (4.12) in (4.15), we obtain

$$\frac{d\kappa_h}{ds_h}\vec{a}_h = 0,$$

which means that $\kappa_h = $ constant. Then from Theorem 3.1 and Theorem 2.1, $S_q$ is an $h$-slant ruled surface in $E^3$.

From Theorem 2.2 and Theorem 4.5 we have the following corollary:

**Corollary 4.1** $S_h$ is a $q$-slant ruled surfaces if and only if $\dfrac{d^2\vec{h}_h}{ds_h^2} + (1+\kappa_h^2)\vec{h}_h = 0$ holds where $\kappa_h$ is the conical curvature of the surface $S_h$.

**Theorem 4.6** Let $S_q$ be a regular ruled surface in $E^3$ with the Frenet frame $\{\vec{q}, \vec{h}, \vec{a}\}$ and with non-zero conical curvature $\kappa_q$. Then, $S_q$ is an $h$-slant ruled surface if and only if the central tangent vector $\vec{a}_h$ of the ruled surface $S_h$ satisfies following differential equation

$$\frac{d^3\vec{a}_h}{ds_h^3} + (1+\kappa_h^2)\frac{d\vec{a}_h}{ds_h} = 0. \qquad (4.16)$$

*Proof.* Let $S_q$ be an $h$-slant ruled surface in $E^3$. From the Frenet formulae in (3.1) we have

$$\frac{d\vec{a}_h}{ds_h} = -\kappa_h \vec{h}_h, \tag{4.17}$$

By taking the derivative of (4.17) we get

$$\frac{d^2\vec{a}_h}{ds_h^2} = -\frac{d\kappa_h}{ds_h}\vec{h}_h - \kappa_h\left(-\vec{q}_h + \kappa_h \vec{a}_h\right) \tag{4.18}$$

Since $S_q$ is an $h$-slant ruled surface in $E^3$, then, $\kappa_h = $ constant, i.e., $\kappa'_h = 0$. Hence, (4.18) turns into

$$\frac{d^2\vec{a}_h}{ds_h^2} = -\kappa_h\left(-\vec{q}_h + \kappa_h \vec{a}_h\right). \tag{4.19}$$

By taking the derivative of (4.19) and using $\kappa'_h = 0$, we obtain

$$\frac{d^3\vec{a}_h}{ds_h^3} + (1+\kappa_h^2)\frac{d\vec{a}_h}{ds_h} = 0,$$

which is desired.

Conversely, let us assume that (4.16) holds. From (3.1) we have

$$\frac{d\vec{a}_h}{ds_h} = -\kappa_h \vec{h}_h. \tag{4.20}$$

By taking the derivative of (4.20) twice and using (4.16) we get

$$\kappa_h \vec{h}_h = 2\frac{d\kappa_h}{ds_h}\vec{q}_h + \left(\kappa_h - \frac{d^2\kappa_h}{ds_h^2}\right)\vec{h}_h - 3\kappa_h \frac{d\kappa_h}{ds_h}\vec{a}_h. \tag{4.21}$$

Since the Frenet frame is linearly independent, from (4.21) we obtain the following system:

$$\frac{d\kappa_h}{ds_h} = 0, \quad \kappa_h \frac{d\kappa_h}{ds_h} = 0, \quad \frac{d^2\kappa_h}{ds_h^2} = 0,$$

which gives us that $\kappa_h = $ constant. Then from Theorem 2.2 and Theorem 3.2, we have that $S_q$ is an $h$-slant ruled surface in $E^3$.

From Theorem 2.1 and Theorem 4.6 we have the following corollary:

**Corollary 4.2** $S_h$ is a $q$-slant ruled surfaces if and only if $\frac{d^3\vec{a}_h}{ds_h^3} + (1+\kappa_h^2)\frac{d\vec{a}_h}{ds_h} = 0$ holds where $\kappa_h$ is the conical curvature of the surface $S_h$.

**Theorem 4.7** *Let $S_q$ be a regular ruled surface in $E^3$ with the Frenet frame $\{\vec{q}, \vec{h}, \vec{a}\}$ and with non-zero conical curvature $\kappa_q$. Then, $S_q$ is a $q$-slant ruled surface if and only if the ruling vector $\vec{q}_a$ of the ruled surface $S_a$ satisfies the following differential equation*

$$\frac{d^3\vec{q}_a}{ds_a^3} + (1+\kappa_a^2)\frac{d\vec{q}_a}{ds_a} = 0. \tag{4.22}$$

***Proof.*** Let $S_q$ be an $q$-slant ruled surface in $E^3$. From the Frenet formulae in (3.5) we have

$$\frac{d\vec{q}_a}{ds_a} = \vec{h}_a.$$

If we take the derivative of the last equation twice with respect to $s_a$, we get

$$\frac{d^3\vec{q}_a}{ds_a^3} = -\frac{d\vec{q}_a}{ds_a} + \frac{d\kappa_a}{ds_a}\vec{a}_a - \kappa_a^2\frac{d\vec{q}_a}{ds_a}. \tag{4.23}$$

Since $S_q$ is a $q$-slant ruled surface in $E^3$, then from Theorem 2.2 and Theorem 3.3, $\kappa_a$ is constant, i.e., $\frac{d\kappa_a}{ds_a} = 0$. Then from (4.23) it follows

$$\frac{d^3\vec{q}_a}{ds_a^3} + (1+\kappa_a^2)\frac{d\vec{q}_a}{ds_a} = 0,$$

which is desired.

Conversely, let (4.22) holds. From (3.5) we have

$$\frac{d\vec{h}_a}{ds_a} = -\vec{q}_a + \kappa_a\vec{a}_a. \tag{4.24}$$

By taking the derivative of (4.24) we get

$$\frac{d^2\vec{h}_a}{ds_a^2} = \frac{d^3\vec{q}_a}{ds_a^3} = -\frac{d\vec{q}_a}{ds_a} - \kappa_a^2\vec{h}_a + \frac{d\kappa_a}{ds_a}\vec{a}_a, \tag{4.25}$$

and using (4.22) in (4.25), we obtain that $\kappa_a = \text{constant}$, which means that $\kappa_q = \text{constant}$ and so $S_q$ is an $q$-slant ruled surface in $E^3$.

**Theorem 4.8** *Let $S_q$ be a regular ruled surface in $E^3$ with the Frenet frame $\{\vec{q}, \vec{h}, \vec{a}\}$ and with non-zero conical curvature $\kappa_q$. Then, $S_q$ is a $q$-slant ruled surface if and only if the central normal vector $\vec{h}_a$ of the ruled surface $S_a$ satisfies the following differential equation*

$$\frac{d^2 \vec{h}_a}{ds_a^2} + (1 + \kappa_a^2) h_a = 0. \tag{4.26}$$

**Proof.** Let $S_q$ be a $q$-slant ruled surface in $E^3$. From the Frenet formulae in (3.5) we have,

$$\frac{d\vec{h}_a}{ds_a} = -\vec{q}_a + \kappa_a \vec{a}_a.$$

By taking the derivative of the last equation we get

$$\frac{d^2 \vec{h}_a}{ds_a^2} = -\frac{d\vec{q}_a}{ds_a} + \frac{d\kappa_a}{ds_a} \vec{a}_a - \kappa_a^2 h_a. \tag{4.27}$$

Since $S_q$ is a $q$-slant ruled surface in $E^3$, from Theorem 2.2, $\kappa_a =$ constant. Then from (4.27) we have

$$\frac{d^2 \vec{h}_a}{ds_a^2} + (1 + \kappa_a^2) h_a = 0.$$

Conversely, let the equation (4.26) holds. From (3.5) we have

$$\vec{q}_a = \kappa_a \vec{a}_a - \frac{d\vec{h}_a}{ds_a}. \tag{4.28}$$

By taking the derivative of (4.28) and using (3.5) and (4.26) we conclude that $\kappa_a =$ constant which means that $S_q$ is a $q$-slant ruled surface in $E^3$.

**Theorem 4.9** Let $S_q$ be a regular ruled surface in $E^3$ with the Frenet frame $\{\vec{q}, \vec{h}, \vec{a}\}$ and with non-zero conical curvature $\kappa_q$. Then, $S_q$ is a $q$-slant ruled surface if and only if the central tangent vector $\vec{a}_a$ of the ruled surface $S_a$ satisfies following differential equation

$$\frac{d^3 \vec{a}_a}{ds_a^3} + (1 + \kappa_a^2) \frac{d\vec{a}_a}{ds_a} = 0. \tag{4.29}$$

**Proof.** Let $S_q$ be a $q$-slant ruled surface in $E^3$. From the Frenet formulae in (3.5) we have

$$\frac{d\vec{a}_a}{ds_a} = -\kappa_a h_a.$$

By taking the derivative of the last equation we get

$$\frac{d^2 \vec{a}_a}{ds_a^2} = -\frac{d\kappa_a}{ds_a} \vec{h}_a + \kappa_a \vec{q}_a - \kappa_a^2 \vec{a}_a. \tag{4.30}$$

Since $S_q$ is a $q$-slant ruled surface in $E^3$, $\kappa_a = $ constant and from (4.30) it follows,

$$\frac{d^2\vec{a}_a}{ds_a^2} = \kappa_a \vec{q}_a - \kappa_a^2 \vec{a}_a. \tag{4.31}$$

By taking the derivative of (4.31) again and using $\frac{d\kappa_a}{ds_a} = 0$ we obtain

$$\frac{d^3\vec{a}_a}{ds_a^3} + (1+\kappa_a^2)\frac{d\vec{a}_a}{ds_a} = 0,$$

which is desired.

Conversely, let us assume that (4.29) holds. From (3.5) we have

$$\frac{d\vec{a}_a}{ds_a} = -\kappa_a \vec{h}_a. \tag{4.32}$$

By taking the derivative of (4.32) twice and using (4.29) we get

$$\kappa_a \vec{h}_a = 2\frac{d\kappa_a}{ds_a}\vec{q}_a + \left(\kappa_a - \frac{d^2\kappa_a}{ds_a^2}\right)\vec{h}_a - 3\kappa_a \frac{d\kappa_a}{ds_a}\vec{a}_a. \tag{4.33}$$

Since the Frenet frame is linearly independent, from (4.33) we obtain the following system:

$$\frac{d\kappa_a}{ds_a} = 0, \quad \kappa_a \frac{d\kappa_a}{ds_a} = 0, \quad \frac{d^2\kappa_a}{ds_a^2} = 0,$$

which leads us to $\kappa_a = $ constant. Then $S_q$ is a $q$-slant ruled surface in $E^3$.